\documentclass[twoside,10pt]{article}

\usepackage{a4}
\usepackage{amsthm}
\usepackage{amsmath}
\usepackage{euscript}
\usepackage{amssymb}

\newtheorem{corollary}{\sc Corollary}[section]
\newtheorem{theorem}{\sc Theorem}[section]
\newtheorem{lemma}{\sc Lemma}[section]
\newtheorem{proposition}{\sc Proposition}[section]
\newtheorem{remark}{\sc Remark}
\newtheorem{definition}{\sc Definition}[section]
\newtheorem{keylemma}{\sc Key-lemma}[section]

\def\R{\mathbb R}
\def\N{\mathbb N}

\def\kac{Ka$\mathrm{\check{c}}$}  

\def\t{{\bf t}}  
\def\r{{\bf r}}  
\def\nor{{\bf{\hat r}}}  
\def\w{{\bf w}}  

\def\A{\EuScript{A}}
\def\renyi{\EuScript{M}}
\def\return{\EuScript{R}}
\def\hreturn{\hat{\EuScript{R}}}

\def\limn{\lim_{n\rightarrow\infty}}           
\def\limsupn{\limsup_{n\rightarrow \infty}}    

\let\goth=\mathfrak
\let\phi=\varphi

\newcommand{\eqdef}{\stackrel{\scriptscriptstyle\rm def}{=}}

\begin{document}

\begin{center}
{\sc{\large Entropy estimation and fluctuations of}}\\
{\sc{\large  Hitting and Recurrence Times}}\\
{\sc{\large for Gibbsian sources}}
\end{center}


\centerline{(Running title: Entropy estimation and fluctuations of hitting times)}

\vskip 1truecm

\centerline{\scshape J.-R. Chazottes$^{a}$, E. Ugalde$^{b}$}
\medskip
 
{\footnotesize
\centerline{$^{a}$ CPHT-CNRS}
\centerline{{\'E}cole Polytechnique, 91128 Palaiseau Cedex, France}
\centerline{{\tt Email: jeanrene@cpht.polytechnique.fr}} 
\centerline{$^{b}$ IICO - UASLP}
\centerline{A. Obreg{\'o}n 64, 78000 San Luis Potos{\'\i}, SLP, M{\'e}xico}
\centerline{{\tt Email: ugalde@cactus.iico.uaslp.mx}}
}

\vskip 2truecm

\begin{abstract}
Motivated by entropy estimation from chaotic time series,
we provide a comprehensive analysis of 
hitting times of cylinder sets in the
setting of Gibbsian sources.
We prove two strong approximation results from which
we easily deduce pointwise convergence to entropy, lognormal
fluctuations, precise large deviation estimates and an explicit
formula for the hitting-time multifractal spectrum.
It follows from our analysis that 
the hitting time of a $n$-cylinder fluctuates in the same way as the
inverse measure of this $n$-cylinder at `small scales', but in a
different way at `large scales'. In particular, the R{\'e}nyi entropy
differs from the hitting-time spectrum, contradicting a naive ansatz.
This phenomenon was recently numerically observed for return
times that are more difficult to handle theoretically.
The results we obtain for return times, though less complete,
improve the available ones.

\bigskip 


\end{abstract}


\newpage

\section{Introduction}

The setting of this work is a model for the following situation: One has time
series $x_1,x_2,x_3,...,x_N,...$ assumed to be typical realizations of some dynamical system.
We further assume that the invariant measure of the system is
Gibbsian. We then want to evaluate the entropy, the R{\'e}nyi
spectrum and the hitting-time/return-time spectrum
of the system by using estimators based on hitting and return times. 
Our aim is to analyse the fluctuation properties of these estimators
to have an {\em a priori} control on what we compute in
practice. Regarding the estimator of the hitting-time/return-time spectrum, the
issue is even to determine to what it converges.

Before being more specific and describing these estimators,
let us temporaly adopt a general point of view on hitting and return times.
The first result in the early Ergodic Theory of dynamical systems is
Poincar{\'e}'s recurrence Theorem, see e.g. \cite{krengel} or virtually any textbook on
ergodic theory.
Colloquially it states that in any dynamical
system preserving a finite measure, typical orbits have a non-trivial
recurrence behavior to each set of positive measure in that they come back
infinitely often to it.
More formally, let $T:X\to X$ be a measurable transformation of the set $X$, the `phase space',
and $\mu$ a $T$-invariant probability measure on $X$. Poincar{\'e}
recurrence Theorem asserts that if $A\subset X$ is a measurable set of
positive measure, then there is a subsequence $\{n_i\}_{i\geq 1}$ of the positive integers
such that $T^{n_i}x\in A$ for $\mu$-almost every point $x\in A$.
One can also ask whether the orbit of a point $x\notin A$ will eventually
enter $A$. A sufficient condition is ergodicity: for $\mu$-almost every $x$, there
is a finite time $n_0=n_0(x)$ such that $T^{n_0}x\in A$. This is the first hitting time.

It is natural to seek for quantitative descriptions of this recurrence
and hitting time.
\kac's Lemma, see e.g. \cite{krengel}, states that the conditional expectation (first
moment) of the first recurrence time to $A$ is equal to $1/\mu(A)$.
This agrees with the intuition that the smaller the measure of $A$ is, the longer it takes
to come back. But ergodicity alone is not sufficient in general to guarantee
the finiteness of the expectation of the first hitting time. Sufficient conditions
on the mixing properties of the dynamics were recently given to ensure the finiteness of finitely many or all
moments of hitting and return times \cite{chazottes1}. We shall use them below.

To go further, we turn our attention to sets $A$ which are cylinder sets.
Let $\A$ be a generating partition of $X$. To each point $x\in X$ we associate
its $n$-cylinder, $\A_n(x)$, defined as the intersection of
all the elements of $\A$, $T^{-1}\A$, ..., $T^{-n+1}\A$ containing $x$.
Let us phrase the following remarkable result \cite{OW}: 
The return time $\r_n(x)$ of a $\mu$-typical point $x$ into  $n$-cylinder $\A_n(x)$,
behaves as follows\footnote{The symbol `$b_n \asymp c_n$' precisely means
$\limn\frac1{n}\log b_n = \limn\frac1{n}\log c_n$ ($\mu$ or
$\mu\!\times\!\mu$ almost surely
for pointwise quantities).}:
\begin{equation}\label{OW}
\r_n(x)\asymp\exp(nh_\mu(T))
\end{equation}
where $h_\mu(T)$ is the measure-theoretic entropy
of the system. An analog result is available for the first time $\w_n(x,y)$ 
the orbit of $x$ enters the $n$-cylinder $\A_n(y)$ about a point $y$, where $x$ and $y$
are `randomly chosen' according to $\mu$ independently of one another:
\begin{equation}\label{shields}
\w_n(x,y)\asymp \exp(nh_\mu(T))
\end{equation}
(But ergodicity is not enough: `strong' mixing properties are necessary \cite{shields}.)
Formulas \eqref{OW} and \eqref{shields} give a simple entropy estimator based
on the observation of a single typical orbit of the system.

It is natural to look for the asymptotic law of return times, that is
to study, for each $t>0$,
\begin{equation}\label{EXPLAW}
\mu\left\{z\in \A_n(x): \r_n(z)\leq \frac{t}{\lambda_{\A_n(x)}}\right\}
\end{equation}
when $n\to \infty$.
An easy heuristic argument shows that one must
rescale $t$ by a quantity proportional to $ \mu(\A_n(x))$
to obtain a non-trivial limiting law.
Typically
\begin{equation}\label{SMB}
\mu(\A_n(x))\asymp \exp(-n h_{\mu}(T))
\end{equation}
by Shannon-McMillan-Breiman Theorem.

Under an ever lessening set of hypotheses on the nature
and the speed of mixing of the system, the limiting law has
been proved to be the exponential law both for hitting and return
times, see for instance \cite{AG,HSV}.  

The purpose of this paper is to analyse the fluctuations of return
times and hitting times by using approximations by the inverse measure of
cylinders.
Indeed, combining \eqref{OW} and \eqref{SMB}, or \eqref{shields} and \eqref{SMB} we get
\begin{equation}\label{rnmu}
\r_n(x)\asymp \frac1{\mu(\A_n(x))}\quad\textup{or}\quad
\w_n(x,y)\asymp \frac1{\mu(\A_n(x))}
\end{equation}
Our aim is to sharpen these rough relations which are only
pointwise. More precisely, we want to compare the fluctuations
of these quantities. Do they have the same log-normal fluctuations ?
The same large deviations ? On another hand, there is a fundamental
class of ergodic measures, namely Bowen-Gibbs measures
\cite{bowen,keller}, such that
\footnote{The symbol $b_n \backsim c_n$ means that $\max(b_n/c_n,
c_n/b_n)$ is bounded from above. We shall use it in the sequel.}
\begin{equation}\label{gibbs}
\mu(\A_n(x)) \backsim \exp(S_n\phi(x))
\end{equation}
where $S_n\phi(x)$ is the `energy' of the cylinder $\A_n(x)$.
Therefore, one can reduce the study of $\log\r_n(x)$
or $\log\w_n(x,y)$ to that of $S_n\phi(x)$, which is much
more easy to tackle and, indeed, all is known on the fluctuations
of $S_n\phi(x)$.

A related issue is to compute the so-called multifractal spectrum or
R{\'e}nyi entropy \cite{P} , defined as
\begin{equation}\label{loosely}
\int \mu(\A_n(x))^{-q}\ d\mu(x) \asymp \exp(nq\overline{\renyi}(q))\quad\textup{as}\;n\to\infty\,.
\end{equation}
In view of practical computation of the R{\'e}nyi entropy, when one only has at hand a time
series, it is tempting to make the ansatz 
$$
\mu(\A_n(x))^{-1}\leftrightarrow\r_n(x)
$$
in \eqref{loosely} and to evaluate the integral as a Birkhoff average.
This was done by Grassberger \cite{G}. The implicit
assumption is that these two quantities have the same {\em large fluctuations}, as we shall explain below.
An even more problematic point is the tacit assumption that {\em all moments} of the Poincar{\'e} recurrences
are {\em finite} before taking the thermodynamic limit. 
The first explicit introduction of such a Poincar{\'e} recurrence spectrum
is done in \cite{HLMV}. On the basis of numerical computations and
heuristic arguments, it is claimed that the Poincar{\'e} recurrence
spectrum and the R{\'e}nyi spectrum do not coincide even in the setting of
Bowen-Gibbs measures : They argue that the Poincar{\'e} recurrence
spectrum must behave like $1/q$ when $q\to-\infty$, which is not the
case for the R{\'e}nyi spectrum. Our goal is to clarify this claim in view
of practical estimation of these spectra.

We will mainly concentrate on hitting times because they are simpler
to analyse and, at the same time, do share the same properties
with return times for strongly mixing measures like Bowen-Gibbs
measures.
Our tools are thermodynamical formalism and 
a very sharp result \cite{miguel} that gives the error term in the convergence
of \eqref{EXPLAW} to the exponential law both in the size of $\A_n$ and in $t$.
From this we derive two approximation results: a global one and a
local one. Indeed, Theorem \ref{SGFWT} gives an approximation, for
{\em any} $n$, of
\begin{equation}\label{approx1}
\int \w_n^q\ d\mu\!\times\!\mu\quad\textup{for all}\;q\in\R
\end{equation}
as certain partition functions; 
Theorem \ref{SA} gives an almost-sure approximation of
\begin{equation}\label{approx2}
\log(\w_n(x,y)\mu(\A_n(x))\,.
\end{equation}
From \eqref{approx1} we deduce our large deviation results and an
{\em explicit formula} for the hitting-time spectrum, whereas from \eqref{approx2} we derive
a central limit theorem and even a law of iterated logarithm.

As a matter of fact, we shall see that $(1/n)\log\w_n(x,y)$ has the
same normal fluctuations as $-(1/n)\log\mu(\A_n(x))$, but {\em their large
deviations are not the same} in some region. This is because
\begin{equation}\nonumber
\int \w_n^{q}\ d\mu\!\times\!\mu \backsim \sum_{\A_n} \mu(\A_n)^{1-q},\quad\textup{for}\;q>-1
\end{equation}
\begin{equation}\nonumber
\int \w_n^{q}\ d\mu\!\times\!\mu \backsim \sum_{\A_n}
\mu(\A_n)^{2},\quad\textup{for}\;q\le -1\,.
\end{equation}
This behavior is numerically observed in \cite{HLMV} for return times.

Concerning return times, we can completely analyse lognormal fluctuations
but not get an explicit formula for the return-time spectrum. Namely we prove that
this spectrum coincides with the R{\'e}nyi spectrum on $[0,+\infty)$.
The reason for this is the presence of `too soon recurrent cylinders'.
Nevertheless, at the end of the paper we study non-overlapping return
times $\nor_n$. For them we prove that $\int\nor_n^q\ d\mu$ becomes
flat for $q<-1$ and coincides with the R{\'e}nyi spectrum for $q\in
[0,\infty[$.
We conjecture that, in fact, the return-time spectrum
really coincides with the hitting-time spectrum for Bowen-Gibbs measures.

\bigskip

{\bf Outline of the paper}.
In Section \ref{setup} we record relevant definitions and results on
hitting times as well as on Bowen-Gibbs measures. 
In Section \ref{MR}, we establish the two main theorems of the paper
for hitting times, namely a strong {\em global} approximation of the `free energy'
of hitting times of $n$-cylinders, and a strong {\em local} approximation of hitting times.
In Section \ref{corollaries}, we derive a number of corollaries from these two theorems:
pointwise convergence, precise large deviation estimates, an explicit formula for the hitting
time spectrum, a central limit theorem and a law of iterated logarithm.
In Section \ref{RT} we deal with return times.
Section \ref{final} contains three subsections. One is about
(non-overlapping) return times. We improve our previous results by considering
non-overlapping return times. Another one is concerned with
bibliographical notes and possible straightforward extensions
of our work. The last one illustrates that for
the Manneville-Pomeau intermittent map, the hitting-time and return-time
spectra are infinite for $q\geq q_c$, where
$q_c>0$ (but we have a finite invariant measure).

\section{Set-up and background}\label{setup}

The phase space $X$ will be the set of sequences $x=(x_1,x_2,...)$ where $x_j\in\mathcal{A}$
(the {\em finite} alphabet), that is $X=\mathcal{A}^{\N}$. The dynamics is given by the shift map $T$
defined as $(Tx)_j=x_{j+1}$ for all $j\geq 1$.
We only consider full shifts since the passage to subshifts
of finite type is not an issue. Given $a_1^n\eqdef a_1 a_2...a_n$, $a_j\in\mathcal{A}$,
we denote by $[a_1^n]$ the corresponding cylinder set that is 
$[a_1^n]=\{x\in X: x_j=a_j,\ j=1,2,...,n\}$. A point $x\in X$ defines a sequence a cylinders
that we naturally denote by $[x_1^n]$, $n\geq 1$. The notation $x_i^j$, $1\leq i\leq j$, stands
for $x_i x_{i+1}...x_j$. The natural $\sigma$-algebra $\goth{B}$ we take is the $\sigma$-algebra generated
by cylinder sets. We omit to mention it in the sequel since it will always be the reference $\sigma$-algebra.

\begin{definition}\label{def-hitting}
We define the (first) hitting time of $x$ to a cylinder $[a_1^n]$ as follows:
$$
\t_{[a_1^n]}(x)\eqdef\inf\{j\geq 1: x_j^{j+n-1}=a_1^n\}\, ,
$$
and the following hitting time, given $x,y\in X$:
$$
\w_n(x,y)\eqdef\t_{[x_1^n]}(y)=\inf\{j\geq 1: y_j^{j+n-1}=x_1^n\}\, 
$$
which is nothing but the (first) time one sees the $n$ first symbols of $x$ appearing in
$y$, i.e., the first time that the orbit of $y$ hits the cylinder $[x_1^n]$.
\end{definition}

The time $\w_n(x,y)$ is also called the waiting time \cite{shields}.

\bigskip

Let us record the useful facts on the class of ergodic measures we are interested in. We refer the
reader to \cite{bowen,pp,walters} for full details.
Let the potential $\phi:X\to\R$ be of `summable variations'. This means
$$
\sum_{n\geq 1} \textup{var}_n\phi <\infty 
$$
where
$\textup{var}_n\phi\eqdef \sup\{|\phi(x)-\phi(y)|: x_1^n=y_1^n, x,y\in X \}$.
The condition imposed in \cite{bowen} is more restrictive since it is
$\textup{var}_n\phi \leq C\theta^n$, for some $C>0$, $0<\theta<1$
(H{\"o}lder continuity).

\noindent {\sc Bowen-Gibbs property}. Assume that $\phi:X\to\R$ has summable variations.
Then there is a unique shift-invariant measure $\mu=\mu_\phi$, that we call
a Bowen-Gibbs measure, such that for all $n\geq 1$, for all $a_1^n$ and for
any $x\in[a_1^n]$
\begin{equation}\label{BGI}
K^{-1}\leq \frac{\mu([a_1^n])}{\exp(-nP(\phi)+S_n\phi(x))}\leq K
\end{equation}
where
$K=K(\phi)>0$, $S_n\phi(x)\eqdef \phi(x)+\phi(T x)+\cdots+\phi(T^{n-1}x)$
and $P(\phi)$ is the topological pressure of $\phi$.
From \eqref{BGI} it is easy to deduce that
$$
P(\phi)=\limn\frac1{n}\log\sum_{a_1^n}\exp(S_n\phi(a_1^n))
$$
where $S_n\varphi(a_1^n)\eqdef \sup\{S_n\varphi(x): x_1^n = a_1^n \}$.
We assume without loss of generality that $\phi$ is normalized,
which implies that$P(\phi)=0$ and $\phi<0$ \cite{walters}. 
A Bowen-Gibbs measure can also be characterized as an equilibrium state, that is 
the (unique) shift-invariant measure $\eta$ that maximizes
$\int \phi\ d\eta + h_\eta (T)$, the maximum being equal to $P(\phi)$.
This is the Variational Principle. 
Since we assume $P(\phi)=0$, this leads to
\begin{equation}\label{normalization}
-\int\phi\ d\mu=h_\mu(T)\ .
\end{equation}

\bigskip

For any $q\in\R$, define
$$
\renyi(q)\eqdef \limn \frac1{n} \log\sum_{a_1^n}\mu([a_1^n])^{1-q}=
\limn \frac1{n} \log\int\mu([x_1^n])^{-q}\ d\mu(x)\,.
$$
Using \eqref{BGI}, we trivially have that
\begin{equation}\label{renyiviaphi}
\renyi(q)=\limn \frac1{n} \log\int e^{-q S_n\phi(x)}\ d\mu(x)\,.
\end{equation}
It can be easily showed that $q\mapsto\renyi(q)$ is a well-defined function on $\R$ 
when $\mu$ is a Bowen-Gibbs measure. Moreover this function is convex and increasing. 
Indeed, by using the Bowen-Gibbs property \eqref{BGI} and the definition of topological
pressure, we easily get $\renyi(q)=P((1-q)\phi)$ (recall that $P(\phi)=0$).

\bigskip

We now state the key result allowing us to analyse fluctuations
of hitting times.
The first result (with its proof) can be found in \cite{miguel}.

\begin{keylemma}[Exponential distribution of hitting times with error term]\label{MKL}
Assume that $\mu$ is a Bowen-Gibbs measure. Then there exist strictly
positive constants $c,C,\rho_1,\rho_2$, with $\rho_1\leq \rho_2$,
such that for all $n\in\N$, all cylinder $[a_1^n]$ and all $t>0$
there exists $\rho(a_1^n)\in[\rho_1,\rho_2]$ such that one has
\begin{equation}\label{strong-approximation}
\Big\vert \mu\left\{z:\t_{[a_1^n]}(z)>\frac{t}{\mu([a_1^n])}\right\}-\exp(-\rho(a_1^n)t)\Big\vert
\leq \varepsilon(a_1^n,t) \;,
\end{equation}
where $\varepsilon(a_1^n,t)\eqdef C \exp(-c n) \exp(-\rho(a_1^n)t)$.
\end{keylemma}

\begin{remark}
The previous result is established under the hypothesis of 
$\psi$-mixing of the process. Bowen-Gibbs measures indeed
have this property.
When the potential is H{\"o}lder continuous, 
the proof of this property is done in \cite{bowen}
and in fact the $\psi$-mixing coefficient decreases exponentially fast.
When the potential has summable but not exponentially small variations, the $\psi$-mixing
property is established implicitly in the proof of Theorem 3.2 in \cite{walters}. Notice
that in this case we do not know how fast the $\psi$-mixing
coefficient decreases.

\end{remark}

\section{Strong approximations of hitting times}\label{MR}

We can state the main theorem of this section. For sequences of real
numbers $(b_n)$, $(c_n)$, the notation
$b_n\backsim c_n$  precisely means that
$\max\{b_n/c_n,c_n/b_n\}$ is bounded from above.

We have the following

\begin{theorem}\label{SGFWT}
Let $\mu$ be a Bowen-Gibbs measure. Then 
\begin{equation}\label{q_large}
\int \w_n^{q}\ d\mu\!\times\!\mu \backsim \sum_{a_1^n} \mu([a_1^n])^{1-q},\quad\textup{for}\;q>-1 
\end{equation}
and
\begin{equation}\label{q_small}
\int \w_n^{q}\ d\mu\!\times\!\mu \backsim \sum_{a_1^n}
\mu([a_1^n])^{2},\quad\textup{for}\;q\le -1 \,.
\end{equation}
\end{theorem}

\begin{proof}
Let $q> 0$. Then
\begin{align}\label{int-repres}
\int \w_n^{q}\ d\mu\!\times\!\mu & =
\sum_{a_1^n} \mu([a_1^n]) \int 
\t_{[a_1^n]}^{q}\ d\mu \\
& = q\sum_{a_1^n} \mu([a_1^n])^{1-q}
\int_{\mu([a_1^n])}^{\infty} t^{q-1} \mu \left\{ \t_{[a_1^n]}>
\frac t{\mu([a_1^n])}\right\} dt.
\end{align} 
By Key-lemma \ref{MKL}
there exist positive constants $A,B$ such
that for any $t>0$ one has
$$
\mu \left\{ \t_{[a_1^n]}> \frac t{\mu([a_1^n])}\right\}\le A e^{-B t}.
$$
Key-lemma \ref{MKL} also easily gives the lower bound~:
$$
\int_{\mu([a_1^n])}^{\infty} t^{q-1} \mu \left\{ \t_{[a_1^n]}>
\frac t{\mu([a_1^n])}\right\} dt\geq
K' - C\exp(-c n)\ K"
$$
where
$0<K'\eqdef \int_1^\infty t^{q-1}\ e^{-\rho_2 t}\ dt<\infty$ and
$0<K"\eqdef \int_0^\infty t^{q-1}\ e^{-\rho_1 t}\ dt<\infty$. 
There exists an integer $n_0$ such that for all $n\geq n_0$, $K' - C\exp(-c n)\ K">0$.
Therefore we get 
$$
K_1\sum_{[a_1^n]} \mu([a_1^n])^{1-q}\le \int \w_n^{q} d\mu\times \mu
\le K_2\sum_{a_1^n} \mu([a_1^n])^{1-q},
$$
where $$
K_1 = q \ (K' - C\exp(-c n_0)\ K")\, ,
 \quad
 K_2 = qA\ \int_{0}^{\infty} t^{q-1} e^{-Bt} dt.
$$
This establishes (\ref{q_large}) for $q\ge 0$. (The case $q=0$ is trivial.)

Let now $q\in (-1,0)$.

\begin{align}
\nonumber
\int \w_n^{-|q|}\ d\mu\times \mu & =
\sum_{a_1^n} \mu([a_1^n]) \int 
\t_{[a_1^n]}^{-|q|}\ d\mu \\
\nonumber
& = \sum_{a_1^n} \mu([a_1^n]) 
\int_{0}^{1} \mu \left\{ \t_{[a_1^n]}^{-|q|}> t\right\} dt\\
\label{integral}
&= |q|\sum_{a_1^n} \mu([a_1^n])^{1+|q|} 
\int_{\mu([a_1^n])}^{\infty} t^{-|q|-1}\ \mu \left\{ \t_{[a_1^n]}\le \frac 
t{\mu([a_1^n])}\right\} dt.
\end{align} 

We first obtain a lower bound for the integral in the last
expression~:
$$
\int_{\mu([a_1^n])}^{\infty} t^{-|q|-1}\ \mu \left\{ \t_{[a_1^n]}\le \frac 
t{\mu([a_1^n])}\right\} dt \geq 
\int_{1}^{\infty} t^{-|q|-1}\ \mu \left\{ \t_{[a_1^n]}\le \frac 
t{\mu([a_1^n])}\right\} dt\,.
$$
Hence, what matters is only the behavior for ``large $t$''. Using
again Key-lemma \ref{MKL}, we get
$$
\mu \left\{ \t_{[a_1^n]}\le \frac 
t{\mu([a_1^n])}\right\} \geq 1- (1+Ce^{-cn})\ e^{-\rho_1 t}\,.
$$
For all $n\geq n_1$, where $n_1\eqdef \frac1{c}
\log\frac{e^{\rho_1}-1}{C}+1$, we have $1- (1+Ce^{-cn})\ e^{-\rho_1}>0$.
Since $1- (1+Ce^{-cn})\ e^{-\rho_1 t} \geq 
1- (1+Ce^{-cn})\ e^{-\rho_1}$ for all $t\geq 1$, we obtain, for all
$n\geq n_1$,
$$
\int_{1}^{\infty} t^{-|q|-1}\ \mu \left\{ \t_{[a_1^n]}\le \frac 
t{\mu([a_1^n])}\right\} dt
\geq 
\frac{(1- (1+Ce^{-cn_1})\ e^{-\rho_1})}{|q|}\,. 
$$
We now turn to the upper bound. We obviously have
$$
\int_{\mu([a_1^n])}^{\infty} t^{-|q|-1}\ \mu \left\{ \t_{[a_1^n]}\le \frac 
t{\mu([a_1^n])}\right\} dt \leq
\left(\int_{0}^{\frac{1}{2}} + \int_{\frac{1}{2}}^{\infty}\right)
 t^{-|q|-1}\ \mu \left\{ \t_{[a_1^n]}\le \frac 
t{\mu([a_1^n])}\right\} dt\,.
$$
The integral from $\frac{1}{2}$ to $\infty$ is finite. 
Now we observe that,
for every $0<t\le\frac1{2}$, we have the following estimate:
$$
\mu \left\{ \t_{[a_1^n]}\le \frac  t{\mu([a_1^n])} \right\}\le
1- e^{-\rho_{2}t}\, .
$$
This estimate follows from the following lemma which is found in
\cite[Lemma 9]{miguel}.
\begin{lemma}\label{pospar}
For any integer  $t$  with $t\mu([a_1^n])\le 1/2$, one has 
$$
\rho_1\le -\frac { \log\mu\{\t_{[a_1^n]}>t\}}{t \mu([a_1^n])}
\le \rho_2\, ,
$$
where $\rho_1$, $\rho_2$ are the (strictly positive) constants
appearing in Key-lemma \ref{MKL}. 
\end{lemma}
Using this estimate for the integral running from $0$ to $\frac{1}{2}$
we get a finite upper bound since 
$$
\int_0^{\frac{1}{2}} \frac{1- e^{-\rho_2 t}}{t^{|q|+1}}\ dt <\infty\,.
$$

Hence we conclude that 
$$
K_1'\sum_{a_1^n} \mu([a_1^n])^{1+|q|}\le \int \w_n^{-|q|} d\mu\times \mu
\le K_2'\sum_{a_1^n} \mu([a_1^n])^{1+|q|},
$$
where $K_1'$ and $K_2'$ are strictly positive constants.
Hence we obtain (\ref{q_large}) for $q\in(-1,0)$.

Finally, let us consider the remaining case $q\le-1$. Then
for sufficiently large $n$ (such that $\mu([a_1^n])<1/2$), one has
\begin{align*}
\int \w_n^{-|q|} d\mu\times \mu &= |q|\sum_{a_1^n} \mu([a_1^n])^{1+|q|}\ 
\int_{\mu([a_1^n])}^{\infty} t^{-|q|-1}\ \mu \left\{ \t_{[a_1^n]}\le \frac 
t{\mu([a_1^n])}\right\} dt\\
&=|q|\sum_{a_1^n} \mu([a_1^n])^{1+|q|}\ 
\Bigl[
\int_{\mu([a_1^n])}^{1/2} 
+\int_{1/2}^{\infty} \Bigr]\
t^{-|q|-1} \mu \left\{ \t_{[a_1^n]}\le \frac 
t{\mu([a_1^n])}\right\} dt\\
&=|q|\sum_{a_1^n} \mu([a_1^n])^{1+|q|} \
\left[I_1(n,a_1^n)  +I_2(n,a_1^n)\right]. 
\end{align*}
Clearly the second integral $I_2(n,a_1^n)$ is uniformly bounded in $n$. Indeed, 
$$
I_2(n,a_1^n)\le \int_{1/2}^\infty \frac 1{t^{1+|q|}}\ dt <+\infty.
$$
However, the first integral $I_1(n,a_1^n)$ is diverging when $n\to\infty$. 
Therefore the limiting behavior as $n\to\infty$ is determined by
$$
|q|\sum_{a_1^n} \mu([a_1^n])^{1+|q|} 
I_1(n,a_1^n) =
$$ 
$$
|q|\sum_{a_1^n} \mu([a_1^n])^{1+|q|} 
\int_{\mu([a_1^n])}^{1/2} 
\mu \left\{ \t_{[a_1^n]}\le \frac 
t{\mu([a_1^n])}\right\}\frac{dt}{t^{1+|q|}}\,.
$$
We again use Lemma \ref{pospar} to get
$$
1- e^{ -\rho_{1}t}\le
\mu \left\{ \t_{[a_1^n]}\le \frac t{\mu([a_1^n])} \right\}
\le 
1- e^{ -\rho_{2}t}\ .
$$
provided that $t\leq \frac1{2}$.
Using \eqref{BGI} (and $P(\phi)=0$) we obtain
$$
K^{-1} \exp(- c'n)\leq \mu([a_1^n])\leq K \exp(-c n)
$$
where $c,c'>0$.
Hence we get
$$
I_1(n,a_1^n)\leq \rho_2 \int_{\mu([a_1^n])}^{\frac1{2}}
t^{-|q|}\ dt
\leq
\frac{\rho_2 (1-2^{|q|-1}K^{-1} e^{-c' n}) }{|q|-1}\ \mu([a_1^n])^{-|q|+1}
$$
where we used the fact that for all $\kappa\in\R$, $1-e^{-\kappa}\le
\kappa$. Notice that for $n$ large enough, the term between parentheses is
strictly positive.
Now, using the fact that $1-e^{-\kappa}\ge \kappa/2$ for any
$\kappa\in[0,1]$, and remembering that $\rho_{1} /2 \le 1$
(\footnote{Indeed, $\rho_2=2$, see \cite{miguel}.}), and using again the Gibbs property
\eqref{BGI}, we obtain
$$
I_1(n,a_1^n)\geq \frac{\rho_1 (1-2^{|q|-1} K e^{-c n})}{2(|q|-1)}\ \mu([a_1^n])^{-|q|+1}
$$ 
where the term between parentheses is strictly positive provided that
$n$ is sufficiently large.
Therefore, for $n$ large enough, we end up with
$$
\frac{|q|\rho_1 (1\!-\!2^{|q|-1}K e^{-c n})}{2(|q|-1)}
\le 
\frac{\int \w_n^{-|q|} d\mu\!\times\!\mu}{\sum_{a_1^n} \mu([a_1^n])^{2}}
\le
\frac{2|q|\rho_2 (1\!-\!2^{|q|-1}K^{-1} e^{-c' n})}{|q|-1}\, .
$$
(Notice that L'H{\^o}pital's rule shows that there is no problem at $q=-1$.)
Thus, we obtain (\ref{q_small}), which finishes the proof.
\end{proof}

We now turn to local strong approximation estimates.

\begin{theorem}\label{SA}
Assume that $\mu$ is a Bowen-Gibbs measure. Then there exists $\epsilon_0>1$ such that for any
$\epsilon>\epsilon_0$, one has
\begin{equation}\label{ae-strong-approximation}
-\epsilon\log n \leq\log\left(\w_n(x,y)\mu([x_1^n])\right)\leq \log\log(n^{\epsilon})
\end{equation}
eventually $\mu\!\times\!\mu$-a.s.
\textup{(}\footnote{By ``eventually $\mu\!\times\!\mu$-a.s.''
we mean that there exists a set $G$ with $\mu\!\times\!\mu(G)=1$ and such that for
any $z\in G$ there is an integer $\tilde{N}=\tilde{N}(z)$ such that for all $n\geq \tilde{N}$
the inequality holds.}\textup{)}.
\end{theorem}

\begin{proof}
We want to find a summable upper-bound to 
$$
\mu\!\times\!\mu\{(x,y):\log(\w_n(x,y) \mu([x_1^n]))> \log t\}=
$$
\begin{equation}\label{split}
\sum_{a_1^n}\mu([a_1^n])\
\mu\left\{\log(\t_{[a_1^n]} \mu([a_1^n]))> \log t \right\}
\end{equation}
where $t$ will be suitably chosen.

First observe that the function $\epsilon(a_1^n,t)\leq C e^{-cn}$ for all $t>0$ in \eqref{strong-approximation}.
Throughout this proof this error bound will be sufficient for our purposes.
Using \eqref{strong-approximation} in \eqref{split} yields
$$
\mu\!\times\!\mu\{(x,y):\log(\w_n(x,y) \mu([x_1^n]))> \log t\}
\leq C e^{-cn}+ e^{-\rho_1 t}\, .
$$
Take $t=t_n=\log(n^{\epsilon})$, where $\epsilon>0$ is to be chosen later on, to get
$$
\mu\!\times\!\mu\{(x,y):\log(\w_n(x,y) \mu([x_1^n]))> \log\log(n^\epsilon)\}
\leq  C e^{-cn} + \frac1{n^{\rho_1 \epsilon}}\,.
$$
Choose $\epsilon>1/\rho_1$.
An application of the classical Borel-Cantelli lemma tells us that
$$
\log(\w_n(x,y) \mu([a_1^n]))\leq \log\log(n^{\epsilon})\quad\textup{eventually a.s.}\,.
$$
For the lower bound, observe that using \eqref{strong-approximation} with the same simplified error bound as
before, we get for all $t>0$
$$
\mu\!\times\!\mu\{(x,y):\log(\w_n(x,y) \mu([x_1^n]))\leq \log t\}
\leq C e^{-cn}+ 1-e^{-\rho_2 t}\leq C e^{-cn}+\rho_2 t\, .
$$
Choose $t=t_n=n^{-\epsilon}$, $\epsilon>1$, to get, proceeding as before,
$$
\log(\w_n(x,y) \mu([x_1^n]))> -\epsilon\log n\quad\textup{eventually a.s.}\,.
$$
The proof is finished by observing that both bounds hold for any $\epsilon>\max(1,\rho_1)$.
\end{proof}

\section{Corollaries}\label{corollaries}

In this section we derive the corollaries of Theorems \ref{SGFWT} and \ref{SA}.

\subsection{Almost-sure convergence}\label{AS}

The following result tells us the `typical' behavior of $\w_n(x,y)$. 
Throughout, $h_\mu (T)$ is the (Kolmogorov-Sinai) entropy of $(X,T,\mu)$.

\begin{corollary}[Almost-sure convergence of hitting times]
Let $\mu$ be a Bowen-Gibbs measure. Then
$$
\limn\frac1{n}\log\w_n(x,y)=h_\mu(T)\quad\textup{for}\,\mu\!\times\!\mu-\textup{a.e.}\,(x,y)\ .
$$
\end{corollary}

\begin{proof}
By Theorem \ref{SA} we get
$$
\limn\frac1{n}\log\w_n(x,y)=\limn-\frac1{n}\log\mu([x_1^n])=
h_\mu(T)\quad\mu\!\times\!\mu-\textup{a.e.}
$$
where the second equality is given by Shannon-McMillan-Breiman Theorem.
\end{proof}

This result means that if we pick up randomly and independently of one another $x$ and $y$, then
the time needed for the orbit of $y$ to hit $[x_1^n]$ is typically of order $\exp(n h_\mu (T))$. 
In fact, this result is valid under the more general assumption that the process is weak Bernoulli
(or $\beta$-mixing). For the details, we refer to \cite{shields}. 
Bowen-Gibbs measures are weak Bernoulli processes, see \cite{bowen} for the proof.

\bigskip

{\bf Remark}. In \cite{chazottes0} the author assumes
that $x$ is picked up randomly according to an ergodic measure $\eta$
whereas $y$ is randomly (and independently)
chosen according to a Bowen-Gibbs measure $\mu$. The previous result becomes:
$$
\limn\frac1{n}\log\w_n(x,y)=h_\eta(T)+h_T(\eta|\mu)\quad\textup{for}\,\eta\!\times\!\mu-\textup{a.e.}\; (x,y)
$$
where $h_T(\eta|\mu)$ is the relative entropy of $\eta$ with respect to $\mu$.
The results obtained in this paper could be suitably generalized to that situation.

\subsection{Large deviations and multifractal spectra}\label{LD}

In this section, we study large deviations of $\frac1{n}\log\w_n$ around the entropy
$h_\mu(T)$ where $\mu$ is a Bowen-Gibbs measure, that is, we only assume that $\phi$ has summable variations.

The ansatz consisting in replacing $\mu([x_1^n])$ in the R{\'e}nyi entropy
by $1/\w_n(x,y)$ leads to the following definition.

\begin{equation}
\mathcal{W}(q)\eqdef\lim_{n\to\infty} \frac{1}{n} \log \int
\w_n^{q}(x,y)\ d(\mu\!\times\!\mu)(x,y)\, ,
\end{equation}
provided the limit exists.
We do not use exactly the same definitions as in \cite{HLMV}. The present definitions
are motivated by large deviation theory.

Introduce, for convenience, the following functions:
$$
\mathcal{W}_n(q)\eqdef \frac1{n} \log \int \w_n^{q}(x,y)\ d(\mu\!\times\!\mu)(x,y)
$$
for all $n\geq 1$ and $q\in\R$. (Notice that $\mathcal{W}_n(q)$ can be infinite.)

Observe that each function $q\mapsto \mathcal{W}_n(q)$
is a convex (increasing) function on $\R$
(hence, in particular, a continuous one
\footnote{Remind that a convex function
defined on a finite interval of the real line can be
discontinuous only at the endpoints of that interval.}).

Observe also that 
$$
\int \w_n^q(x,y)\ d(\mu\!\times\!\mu)(x,y) =  \sum_{a_1^n}\mu([a_1^n])\int \t_{[a_1^n]}^q(y)\ d\mu(y)\\
$$

There is no \kac formula for hitting times for general ergodic
dynamical systems. Ergodicity only ensures that almost surely there is
a finite first hitting-time.
The only fact we know without any assumption is that
$\int \t_{[a_1^n]}^{q}\ d\mu<\infty$ for all $q\leq 0$. Indeed, 
$$
\int \t_{[a_1^n]}^{q}\ d\mu\leq 1
$$
for any $q\leq 0$.
From \cite{chazottes1} it follows that for
any Bowen-Gibbs measure and for any cylinder $[a_1^n]$, we have
$$
\int \t_{[a_1^n]}^{q}\ d\mu<\infty
\quad\textup{for all}\;q\in\R\, .
$$
Hence $\mathcal{W}_n(q)<\infty$ for all $q\in\R$, $n\geq 1$.

We now turn to large deviation results. We refer the reader to \cite{DZ} for background on this
topic.

\begin{corollary}[Scaled generating function of hitting times]
Assume that $\mu$ is a Bowen-Gibbs measure. Then
\begin{equation}\label{sgf-of-w}
\mathcal{W} (q)= \begin{cases} 
\renyi(q), & \textup{ for }q\geq -1,\\
P(2\varphi), &\textup{ for }q<-1,
\end{cases}
\end{equation}
If $\mu$ is not the measure of maximal entropy, then the function
$q\mapsto\mathcal{W}(q)$ is strictly convex on $(-1,\infty)$.
\end{corollary}

\begin{proof}
Clearly (\ref{q_large}) and (\ref{q_small}) imply (\ref{sgf-of-w}).
\end{proof}

Notice that $q\mapsto \mathcal{W}(q)$ is continuous (as it must be for a convex function on $\R$)
but not differentiable at $q=-1$. Indeed, it can be easily checked
that the right derivative at $-1$ of $\mathcal{W}$ is not zero but equal to $-\int\phi\
d\mu_{2\phi}>0$, where $\mu_{2\phi}$ is the Bowen-Gibbs measure for the
potential $2\phi$.

\bigskip

\begin{corollary}[Large deviations of $\w_n$]\label{LDWT}
Let $\mu$ be a Bowen-Gibbs measure which is not the measure of maximal entropy.
Then for all $u\ge 0$ we have
$$
\limn\frac1{n}\log(\mu\!\times\!\mu)\left\{\frac1{n}\log\w_n >h_\mu(T)+u\right\}= 
\inf_{q>-1} \left\{-(h_{\mu}(T)+u) q + \mathcal{W} (q)\right\}
$$
and for all $u\in (0, u_0)$, $u_0\eqdef |\lim_{q\downarrow -1} \mathcal W'(q)-h_{\mu}(T)|$,
$$
\limn\frac1{n}\log(\mu\!\times\!\mu)\left\{\frac1{n}\log\w_n <h_\mu(T)-u\right\}=
\inf_{q>-1} \left\{-(h_{\mu}(T)-u) q + \mathcal{W} (q)\right\}
$$
\end{corollary}

Notice that $u_0> h_\mu(T)$, that is, we capture the large fluctuations of
$\log\w_n/n$ above and below $h_\mu(T)$ since
$\mathcal W'(0)=h_{\mu}(T)$ (see the appendix for the proof).

\bigskip

\begin{proof}
By Theorem \ref{SGFWT} and \eqref{renyiviaphi} we immediately get that for any $q>-1$
$$
\mathcal{W}(q)=\renyi(q)=\limn\frac1{n}\log\int\exp(-qS_n\phi)\ d\mu\,.
$$
It can be easily deduced from \cite{TV} that the function $q\mapsto P(q\varphi)$ is $C^1$
and strictly convex if (and only if) $\mu$ is not the measure of maximal entropy.
In the H{\"o}lder continuous case, it is real analytic and also strictly convex
if (and only if) $\mu$ is not the measure of maximal entropy \cite{pp}.

We can apply a large deviation result due to Plachky and Steinebach \cite{PS}.
(Recall that a strictly convex differentiable function has a strictly increasing
derivative.)
\end{proof}

Let us remark that when the measure is the one of maximum entropy, there are no large fluctuations which is
not surprising.

The R{\'e}nyi spectrum is defined here as $\overline{\renyi}(q)\eqdef\renyi(q)/q$
and the hitting-time spectrum
$\overline{\mathcal{W}}(q)\eqdef \mathcal{W}(q)/q$.
We get $h_\mu(T)$ for the value of these spectra at $q=0$ (using L'H{\^o}pital's rule).

\begin{corollary}\label{multifractal-spectra}
For any Bowen-Gibbs measure we have the following:
$$
\overline{\mathcal{W}}(q)=\overline{\renyi}(q)\quad\textup{for}\;q\geq -1\quad\textup{and}\quad
\overline{\mathcal{W}}(q)=P(2\phi)/q\quad\textup{for}\;q<-1\,.
$$
\end{corollary}

\subsection{Log-normal fluctuations}\label{lognormal}

The purpose of this section is to show that $\w_n(x,y)$ and $1/\mu([x_1^n])$
have the same lognormal fluctuations for Bowen-Gibbs
measures associated to H{\"o}lder continuous potentials.
Namely, we prove a central limit theorem and a law of iterated
logarithm.

We refer the reader to \cite{pp} for full details on the central limit theorem for
H{\"o}lder continuous observables with respect to Bowen-Gibbs measures
with a H{\"o}lder continuous potential.
Define the following variance:
\begin{equation}\label{defvariance}
\sigma^2  \eqdef \limn\frac1{n}\int \left(S_n\phi - n\int\phi\ d\mu\right)^2\ d\mu\,.
\end{equation}
It is well-known that if $\sigma^2>0$ ($\sigma^2<\infty$ because of the exponential decay of correlations) and
\begin{equation}\label{phiclt}
\forall t\in\R\quad\limn\mu\left\{\frac{-S_n\phi -n h_\mu(T)}{\sigma\sqrt{n}}<t\right\}=
\mathcal{N}(0,1)((-\infty,t])
\end{equation}
where
$\mathcal{N}(0,1)((-\infty,t])=\frac1{\sqrt{2\pi}} \int_{-\infty}^{t} \exp(-\xi^2 /2)\ d\xi$.
This is just a particular instance of the central limit theorem for Bowen-Gibbs measures \cite{pp}
where the observable is $-\phi$. By \eqref{normalization} $\int -\phi\ d\mu=h_\mu(T)$.
One has $\sigma^2=0$ if and only if $-\phi+h_\mu(T)$ (or equivalently
$\phi-\int\phi\ d\mu_\phi$) is a coboundary,
i.e. a function of the form $\varrho-\varrho\circ T$,
for some measurable function $\varrho$. 
This means that $\sigma^2=0$ if and only if
$\mu$ is the (unique) measure of maximal entropy.

\begin{corollary}\label{CLTW}
Assume that $\mu$ is a Bowen-Gibbs measure with a H{\"o}lder continuous potential
which is not the measure of maximal entropy. Then
\begin{equation}
\limn\mu\!\times\!\mu\left\{\frac{\log \w_n -n h_{\mu}(T)}{\sigma\sqrt{n}}<t\right\}=
\mathcal{N}(0,1)((-\infty,t])\, .
\end{equation}
Moreover,
\begin{equation}\label{w-variance}
\sigma^2=\limn\frac1{n} \int (\log\w_n - h_\mu(T))^2\ d(\mu\!\times\!\mu)\,.
\end{equation}
\end{corollary}
\begin{proof}
Our goal is to show that the central limit theorem for $\log\w_n(x,y)$ results from
the one for $-\log\mu([x_1^n])$ which in turn results from the one for $-S_n \phi(x)$.
The latter assertion is trivial since, for all $x$,
\begin{equation}\label{Amuphi}
-S_n \phi(x)- C \leq -\log\mu([x_1^n])\leq -S_n \phi(x)+ C
\end{equation}
where $C>0$, due to Bowen-Gibbs inequality \eqref{BGI}.
Hence the quantities $-\log\mu([x_1^n])$ and $-S_n \phi(x)$
have the same mean $h_\mu(T)$ and variance $\sigma^2$. 

We now use the strong approximation formula \eqref{ae-strong-approximation} from Theorem \ref{SA} together
with inequalities \eqref{Amuphi} to get
\begin{equation}\label{AEphir}
\frac{\log\w_n(x,y)+S_n \phi(x)}{\sigma\sqrt{n}}\to 0\quad\textup{for}\,
\mu\!\times\!\mu-\textup{almost every}\, (x,y)\, .
\end{equation}
By a basic result of Probability Theory (see \cite{durrett} for instance), the preceding result
and \eqref{phiclt} imply the desired statement.

The proof of \eqref{w-variance} is given in the Appendix. 
\end{proof}

We emphasize that \eqref{w-variance} does not follow from Theorem
\ref{SA}, see the Appendix.

We now state and prove a law of iterated logarithm for $\log\w_n$:

\begin{corollary}\label{lil}
Assume that $\mu$ is a Bowen-Gibbs measure with a H{\"o}lder continuous
potential which is not the measure of maximal entropy. Then
\begin{equation}\label{lilforw}
\limsupn \frac{\log\w_n - n h_\mu(T)}{\sigma\sqrt{2n\log\log n}}=1\quad\mu\!\times\!\mu-\textup{a.e.}
\end{equation}
\end{corollary}

Remark that we get $-1$ instead of $1$ when taking `$\liminf$' instead
of `$\limsup$'. In fact, we could show that the set of accumulation
points of the sequence $\{(\log\w_n-n h_\mu(T))/\sqrt{2n\log\log
  n}\}_n$ is the interval $[-\sigma,+\sigma]$. 

\bigskip

\begin{proof}
Using Theorem \ref{SA} we get that eventually $\mu\!\times\!\mu$-almost surely
$$
-\frac{\epsilon\log n}{\sigma\sqrt{2n\log\log n}}\leq
\frac{\log\w_n - n h_\mu(T)}{\sigma\sqrt{2n\log\log n}}-
\frac{-S_n\phi - n h_\mu(T)}{\sigma\sqrt{2n\log\log n}}
\leq \frac{\log\log n^\epsilon}{\sigma\sqrt{2n\log\log n}}
$$
Taking the limit supremum $n\to\infty$ and using the law of iterated logarithm for $-S_n\phi$,
we finish the proof.
\end{proof}

The law of iterated logarithm for $-S_n\phi$ can be found in \cite{DP}. 

\section{Return times}\label{RT}

We now turn to return times. As we shall see, we obtain less complete
results than for hitting times.

\subsection{Set-up}

\begin{definition}\label{def-basicreturn}
The (first) return time of a point $x$ into its $n$-cylinder $[x_1^n]$, $n\geq 1$ is defined
as:
$$
\r_n(x)\eqdef\inf\{k\geq 2: x_k^{k+n-1}=x_1^n\}\, .
$$
\end{definition}

The following result is proved in \cite[Section 6]{miguelnew}. In order to state
it we need to define the set of $n$-cylinders with `internal periodicity' $p\leq n$:
$$
\mathcal{S}_p(n)\eqdef \{[a_1^n]: \min\{k\in\{1,...,n\}: [a_1^n]\cap T^{-k}[a_1^n]\neq\emptyset\}=p \}\,.
$$
Notice that the set of $n$-cylinders can be written as the union $\bigcup_{1\leq p\leq n} \mathcal{S}_p(n)$.

\begin{keylemma}[Exponential distribution of basic return times]\label{KL2}
Let $\mu$ be a Bowen-Gibbs measure. Then there exist strictly
positive constants $\tilde{C}_1,\tilde{C}_2,\tilde{C}_3$ 
such that for any $n\in\N$, any $p\in\{1,...,n\}$,\
any cylinder $[a_1^n] \in \mathcal{S}_p(n)$, one has for all $t\geq p$
\begin{equation}\label{strong-approximation-for-basicreturns}
\Big\vert \mu
\left\{z:\r_{n}(z)>\frac{t}{\zeta(a_1^n)\mu([a_1^n])}\ \big\vert\ [a_1^n]\right\}- \zeta(a_1^n) \exp(-t)\Big\vert
\leq \tilde{\varepsilon}(a_1^n,t) \;,
\end{equation}
where $\tilde{\varepsilon}(a_1^n,t)\eqdef \tilde{C}_1 \ \exp(- \tilde{C}_2 n) \ \exp(-\tilde{C}_3 t)$
and $\zeta(a_1^n)$ is such that $|\ \zeta(a_1^n)-\rho(a_1^n)|\leq D\ \exp(- \tilde{C}_2 n)$, with
$D>0$ and $\rho(a_1^n)$ given in Key-lemma \ref{MKL}.
Moreover,
$$
\mu\{z:\r_n(z)>t\ | \ [a_1^n]\}=1\quad \textup{for all}\quad t<p\, .
$$
\end{keylemma}

\subsection{Large deviations}

We define the following (possibly infinite) quantities, for all $\in\R$, provided the limit exists:
\begin{equation}
\return(q)\eqdef \limn \frac1{n} \log \int \r_n^{q}(z)\ d\mu(z)
\end{equation}
$$
\return_n(q)\eqdef \frac1{n} \log \int \r_n^{q}(z)\ d\mu(z)\;\textup{and}
$$
Without any assumption on $\mu$, the function $\return$ trivially exists at $q=0$ 
and equal zero.
If $\mu$ is assumed to be ergodic then $\return(1)=\log|\mathcal{A}|$.
Indeed,
$$
\int \r_n(z)\ d\mu(z)=
\sum_{a_1^n}\left(\int\r_n(z)\ d\mu_{[a_1^n]}(z)\right)\ \mu([a_1^n])\, .
$$
By \kac's formula
the integral is equal to $1/\mu([a_1^n])$, hence we get
$$
\limn\frac1{n}\log\int \r_n(z)\ d\mu(z)=\limn\frac1{n}\log\sum_{a_1^n} 1 =\log|\mathcal{A}|.
$$

Let us emphasize that the finiteness of $\return$ 
for all $q>1$ is not obvious at all.
We need to know the finiteness of $\int_{[a_1^n]}\r_n^{q}\ d\mu$ for all $q>1$, that is the
finiteness of the moments of return times  to $[a_1^n]$.
This point seems to have been overlooked before.
For $q\leq 1$, all moments of return times are of course finite due to \kac's formula
($T$-invariance is indeed sufficient)
but nothing rules out {\it a priori} the possibility that the moment of the return time be infinite beyond
a certain order $q>1$ for some $n_0$ (and hence for all $n\geq n_0$ since the moment of order $q>0$ as a
function of $n$ is increasing). This will be illustrated in Section \ref{MP}.
From \cite{chazottes1} it follows that
for any Bowen-Gibbs measure and for any cylinder $[a_1^n]$, we have
$$
\int_{[a_1^n]} \r_n^{q}\ d\mu_{[a_1^n]}<\infty\quad\textup{for all}\;q\in\R\, .
$$
Hence $\return_n(q)<\infty$ for all $q\in\R,n\in\N$. 

\begin{proposition}[Partial large deviations for $\r_n$]\label{LDbasicRT}
Let $\mu$ be a Bowen-Gibbs measure which is not the measure of maximal entropy.
Then 
$$
\limn\frac1{n}\log\mu\left\{\frac1{n}\log\r_n >h_\mu(T)+u\right\}= 
\inf_{q\geq 0} \left\{-(h_{\mu}(T)+u) q + \return(q)\right\}
$$
that holds for all $u\ge 0$.
\end{proposition}

\begin{proof}
From Key-Lemma \ref{KL2}, we can deduce that
\begin{equation}\label{pouet}
\int \r_n^{q}\ d\mu \backsim \sum_{a_1^n} \mu([a_1^n])^{1-q},
\end{equation}
for all $q>0$ (the case $q=0$ is trivial). This implies that
$\renyi(q)=\return(q)$ for all $q\geq 0$. Thus we can again apply
Plachky-Steinebach's large deviation estimate on $\R^+$.
We leave the details of the proof of \eqref{pouet} to the reader since
it is very similar to the one for hitting times we gave above.
\end{proof}

Some large deviation results are given in \cite{CGS} for $\r_n$,
but they are valid only in a small (non-explicit)
interval around $h_\mu(T)$. This restriction is due to `too soon
recurrent cylinders'. Because of the same problem, we can only extend
this result to the whole range of possible deviations above $h_\mu(T)$.

\subsection{Lognormal fluctuations}

We summarize what happens for return times in the following theorem and corollary.

\begin{theorem}\label{SA-for-rn}
Assume that $\mu$ is a Bowen-Gibbs measure. Then there exists $\epsilon_1>1$ such that for any
$\epsilon>\epsilon_1$, one has
\begin{equation}
-\epsilon\log n \leq\log\left(\r_n(x)\mu([x_1^n])\right)\leq \log\log(n^{\epsilon})
\end{equation}
eventually $\mu$-a.s.
\end{theorem}

\begin{corollary}\label{lognormalRT}
Let $\mu$ be Bowen-Gibbs measure
associated to H{\"o}lder continuous potential.
Then Corollaries \ref{CLTW}-\ref{lil} hold for
$\r_n$ instead of $\w_n$.
($\mu\!\times\!\mu$ is replaced by $\mu$.)
\end{corollary}

{\bf Sketch of proofs}. The proof of the corollary follows exactly the same lines as
that for hitting times. 
Let us sketch the upper bound, leaving the lower bound to the reader.
To apply the classical Borel-Cantelli lemma, as before, we need to upper bound
$$
\mu\left\{x:\log(\r_n(x) \mu([x_1^n]))\geq \log t\ | \ [a_1^n]\right\}=
$$
\begin{equation}\label{UBforBR}
\sum_{p=1}^{n}\quad \sum_{[a_1^n]\in \mathcal{S}_p(n)}\mu([a_1^n])\
\mu\big\{\log(\r_n \mu([a_1^n]))\geq \log t \ | \ [a_1^n]\big\}
\end{equation}
where $t$ will be suitably chosen.
We now use Key-Lemma \ref{KL2} to get
$$
\eqref{UBforBR} \leq 
\sum_{p=1}^{n}\quad \sum_{[a_1^n]\in \mathcal{S}_p(n)}\mu([a_1^n])\
\left( (\rho_2+D) e^{-\frac{\rho_1}{2}t} + \tilde{C}_1 e^{-\tilde{C}_2 n}
\right)
$$
for all $t\geq p\ \mu([a_1^n])$, where we used the fact that there is some $n_0$ such 
that for all $n\geq n_0$,
$\rho_1/2 \leq \rho_1 - D e^{-\tilde{C}_2 n} \leq \zeta(a_1^n)\leq \rho_2 + D$.
Now choose $t=t_n=\log(n^\epsilon)$, $\epsilon>0$
and notice that $t_n\geq p\ \mu([a_1^n])$, for all $p=1,...,n$ for $n$ large enough
since $\mu([a_1^n])\leq K\ e^{-cn}$ by the Bowen-Gibbs property.
Therefore we get
$$
\eqref{UBforBR} \leq 
\underbrace{\sum_{p=1}^{n}\quad \sum_{[a_1^n]\in \mathcal{S}_p(n)}\mu([a_1^n])}_{=1}\
\left( \frac{\rho_2+D}{n^{\epsilon\rho_1 /2}} + \tilde{C}_1 e^{-\tilde{C}_2 n}
\right)
$$
which is summable in $n$ provided that $\epsilon>2/\rho_1$.
We leave the rest of the proof to the reader.

\section{Final comments and open questions}\label{final}


\subsection{Hitting-time and return-time spectra in presence of intermittency}\label{MP}

We consider here the following intermittent map (the so-called
Manneville-Pomeau map) defined on the interval $[0,1]$:
$$
T:x\mapsto x + x^{1+\alpha}\quad\textup{mod}\; 1,\quad\alpha\in(0,1)\;.
$$
The thermodynamic formalism for such a map is now well-understood. We refer the reader
to the recent paper \cite{frank} for more details and references on what we will use.

Let $\phi=-\log|T'|$, the potential function. The map $T$ admits an absolutely continuous
invariant measure $\mu$ which is an equilibrium state for the potential $\phi$.
($\mu$ is not the only equilibrium state. Any measure of the form $t\mu+(1-t)\delta_0$,
where $t\in[0,1]$ is an equilibrium state for $\phi$; $\delta_0$ is the Dirac measure
at $0$.)
This means that $P(\phi)=h_\mu (T)+\int \phi\ d\mu$, where $P(\phi)$ is the topological pressure
of $\phi$ and $h_\mu(T)$ is the measure-theoretic entropy. In fact $P(\phi)=0$ because of the
Rokhlin formula. It is well-known that $\renyi(q)=P((1-q)\phi)$. From the behavior of the pressure
function we get the following properties for $\renyi(q)$:
it is continuous, convex and non-decreasing.
Moreover, $P((1-q)\phi)=0$ for $q\leq 0$, $P((1-q)\phi)>0$ for $q>0$ and $q\mapsto P((1-q)\phi)$
is a real-analytic function for $q>0$. At the critical point one has the following
asymptotics:
$$
\frac{P((1-q)\phi)}{q}\to h_\mu(T)\;\textup{as}\;q\downarrow 0\;.
$$

The Manneville-Pomeau map has two intervals of monotonicity, $I_0$, $I_1$, from which one can define cylinder sets:
$I_{i_1,i_2,...,i_n}(x)=I_{i_1}\cap T^{-1}I_{i_2}\cap\cdots\cap T^{-1}I_{i_n}$ is that interval
of monotonicity for $T^n$ which contains $x$. Contrarily to the case when maps are everywhere expanding,
the ratio
\begin{equation}\label{weak-BGI}
\frac{\mu(I_{i_1,i_2,...,i_n}(x))}{\exp\left( \sum_{k=0}^{n-1}\phi(T^k x)\right)}
\end{equation}
is not uniformly bounded in $n$ and $x$. This comes from the fact that distortions are not
bounded. A more careful analysis shows that one can find bounds from above and below which
are polynomial in $n$ and uniform in $x$. Such a measure is an example of a
weak Gibbs measure.

The following basic proposition shows that large deviation results in
the usual sense do not hold.

\begin{proposition}\label{MP-results}
For all $q\geq\frac1{\alpha}$, we have for all $n\geq 1$
that $\mathcal{W}_n(q)=\infty$ (hence $\mathcal{W}(q)=\infty$,
but $\renyi(q)$ is finite for every $q\in\R$. The same occurs for $\return_n(q)$
for all $q\geq \frac1{\alpha}+1$.
\end{proposition}

\begin{proof}
We are going to show that $\mathcal{W}_1(q)=\log\int \w_1^q\ d(\mu\!\times\!\mu)$ becomes infinite
from some $q_0=q_0(\alpha)>0$ on (hence $\mathcal{W}_n(q)=\infty$ for all $n$ when $q\geq q_0$ since 
$\w_{n+1}(x,y)\geq \w_n(x,y)$ for all $(x,y)$.)
The main point is the following estimate for $\mu(I_{00...0})$:
$$\mu(I\!\!_{\underbrace{0...0}_{\ell\,symbols}})\geq C\ \ell^{-\frac{1}{\alpha}}
$$
where $C$ is a positive constant. But 
$\mu(I_{0\cdots 0})$ (with $\ell$ symbols)
is nothing but the measure of points that do not enter the right interval $I_1$ before 
$\ell$ iterations.
Now observe that
$$
\int \w_{1}^q\ d\mu\!\times\!\mu\geq \ \mu(I_1)\ \int \t_{I_1}^q\ d\mu=
\mu(I_1)\ \sum_{\ell\geq 0}[(\ell+1)^{q}-\ell^q]\
\mu(I\!\!_{\underbrace{0...0}_{\ell\,symbols}})\; .
$$
Therefore $\mathcal{W}_n(q)=\infty$ for all $n\geq 1$ and $q\geq\frac1{\alpha}$.

To pass from $\mathcal{W}_n(q)$ to
$\return_n(q)$ use Proposition 1 in \cite{chazottes1}.
\end{proof}

It was recently proved in \cite{FHV} that, for $0<\alpha<\frac{1}{2}$,
$1/\mu(I_{i_1,i_2,...,i_n}(x))$ and $\r_n(x)$ have the same lognormal fluctuations

Proposition \ref{MP-results} shows that the return and hitting time
spectra are not relevant for non-uniformly hyperbolic systems.
Indeed, a single indifferent fixed point makes these spectra infinite for all $q\geq q_c(\alpha)$.
Moreover, it implies that $1/\mu(I_{i_1,i_2,...,i_n}(x))$ and
and $\r_n(x)$ cannot have the same large deviations.

\subsection{More on return times}

In \cite{HLMV} the authors study the recurrence-time spectrum $\return(q)$.
They show heuristically that $\return(q)$ must behave like $1/q$ as $q\to\infty$.
(Be careful of the different convention used therein to define $\return(q)$.)
Two numerical simulations confirm this behavior: 
The graphs of the spectrum really look like a constant divided by $q$
for $q\gtrsim 1$.
This is indeed what we get rigorously (remember that $q=2$ corresponds to $q=-1$ in our convention)
for hitting-time spectrum $\mathcal{W}(q)$. We are not able, as we saw above, to prove this for return times. 
The difficulty comes from `too soon recurrent' cylinders. More technically speaking, we do not
have the analog of Lemma \ref{pospar} for return times and therefore the finiteness
of the integral in \eqref{integral} escapes us. 

We can go a bit further by looking at non-overlapping return times $\nor_n$ (studied
recently in \cite{maurer}). By definition,
such return times cannot be `too small':
The (first) non-overlapping return time of a point $x$ into its $n$-cylinder
$[x_1^n]$, $n\geq 1$ is defined as:
$$
\nor_n(x)\eqdef\inf\{k\geq 1: x_{kn +1}^{(k+1)n}=x_{1}^{n}\}\, .
$$

We have the following approximation result:
\begin{proposition}
Let $\mu$ be a Bowen-Gibbs measure. Then, for every $q<-1$,
$$
\int \nor_n^q\ d\mu\backsim  \sum_{a_1^n} \mu([a_1^n])^2\,.
$$
For every $q\geq 0$, 
$$
\int \nor_n^q\ d\mu\backsim  \sum_{a_1^n} \mu([a_1^n])^{1-q}\,.
$$
\end{proposition}
(The symbol $\backsim$ is precisely defined at the beginning of Section \ref{MR}.)
Therefore, if we let $\hreturn(q)$ be the analog of $\return(q)$ 
where $\nor_n$ replaces $\r_n$, the previous result implies
that
$$
\hreturn(q)=P(2\phi)\quad\textup{for all}\;q<-1\quad\;\textup{and}\;
\hreturn(q)=\renyi(q)\quad\textup{for all}\;q\geq 0
$$
(Use the Bowen-Gibbs property and the definition of topological pressure
of Section \ref{setup}.) 

Recall that for hitting times we proved a  more precise result (Theorem \ref{SGFWT})
since the second approximation works for $q\in [-1,\infty[$ in that case.

\bigskip

\begin{proof}
First write
\begin{equation}\label{decomp}
\int \nor_n^q\ d\mu=\sum_{k=1}^{\infty} k^{-|q|}\ \mu\{\nor_n = k\}\,.
\end{equation}
Using the $\psi$-mixing property (see \cite{maurer} for details) we get
that 
$$
\mu([a_1^n])(1-\psi(k-1))\leq\mu_{[a_1^n]}\{\nor_n=k\}\leq \mu([a_1^n])(1+\psi(k-1))\,. 
$$
Hence,
$$
\frac{1}{2}\sum_{a_1^n} \mu([a_1^n])^2\leq
\mu\{\nor_n=k\}
\leq \frac{3}{2}\sum_{a_1^n} \mu([a_1^n])^2\,.
$$
Therefore we get
$$
\frac1{2} C_q
\leq\frac{\int \nor_n^q\ d\mu}{\sum_{a_1^n} \mu([a_1^n])^2}\leq
\frac{3}{2} C_q
$$
where $C_q\eqdef \sum_{k=1}^\infty k^{-|q|}<\infty$
since $|q|>1$.
We leave the proof of the statement in the range $q\in[0,\infty[$ to the reader
(use the analog of Key-lemma \ref{KL2} for $\nor_n$ that can be extracted from
\cite{maurer}).

So, we arrived at the desired result.
\end{proof}

We conjecture that Theorem \ref{SGFWT} is true with return-times instead of hitting-times.
It could be easier to first prove this conjecture for non-overlapping return times.

\subsection{Relevance of the hitting-time and return-time spectra}

In view of \cite{feng} (saturation of level sets)
and what we obtained in the present paper, one
can legitimately ask what is the relevance of the return
time spectrum,
except as a trick to compute the R{\'e}nyi spectrum for $q>-1$.
The same can be said for the hitting-time spectrum.
Even in the comfortable setting of Bowen-Gibbs measures,
these spectra contain no information for $q\leq -1$.
In presence of intermittency, we saw that they are
infinite for $q\geq q_c$. 

\subsection{Related works and an extension}

It is worth to indicate to the reader the differences between our work
and the previous ones. 
In the paper \cite{CGS}, the authors study only
return times for Bowen-Gibbs measures associated to a H{\"o}lder continuous potential.
They prove a central limit theorem and a large deviation principle. 
Here we not only improve the lognormal approximation but also extend the range
of accessible large deviations above the true entropy. Moreover we handle potentials with summable
variations and not only H{\"o}lder continuous ones for large deviations.
We also mention \cite{saussol} for a general (but much less sharp) result 
on lognormal fluctuations for return times. 
In the context of $\psi$-mixing stochastic processes, there are two references \cite{konto,wyner}.
Both deal with local strong approximations, in particular central limit theorems. 
The author of \cite{konto} directly uses the $\psi$-mixing property. 
The author of \cite{wyner} first proves an approximation to the
exponential law of rescaled hitting and return times and then deduces 
strong local approximations. We emphasize that his approximation is much less sharp
than the one we use here. This difference is not relevant for deriving strong local approximations
but becomes essential to handle large deviations.
The only paper dealing with large deviations of hitting times is \cite{ACRV} where the authors
study the first occurrence of a cylindrical pattern in the realisation of a Gibbsian random
field on the lattice $\mathbb{Z}^d$, $d\geq 2$. Our proof is very similar to that of this work.

We also note that our results can be extended to a more general class of processes, namely
the processes satisfying the $\phi$-mixing property with a summable $\phi$-mixing sequence.
This is because Key-lemma \ref{MKL} is proved not only for $\psi$-mixing but also
for such processes \cite{miguel}. But up to our knowledge, this does not define a
natural class of equilibrium states on shift spaces. That is why we did not state our results under
this assumption. Nevertheless, an interesting class of non-Markov expanding maps of
the interval has this property with an exponentially decreasing $\phi$-mixing sequence
(with respect to the partition given by the discontinuity points of the map).
This class was indeed studied in \cite{paccaut}. We could therefore sharpen the results
of that paper since Corollary \ref{lognormalRT},
apply. We could of course write down the analogous results for hitting times.
Regarding large deviations of hitting times, we could derive some approximations in the spirit of
Theorem \ref{SGFWT} and derive some estimates like that of Corollary \ref{LDWT}. But one has
to be careful with the control of some `bad' cylinders for which the "distorsion property"
(the analog of Bowen-Gibbs property \eqref{BGI} in that context) does not hold, which is
the price to be paid for the non-Markovian partition.

\section{Appendix}\label{L1L2}

We prove the convergence in mean and in quadratic mean 
of $\log\w_n$.
The former is related to the
slope of $\mathcal{W}$ at $0$, see the comment
after Corollary \ref{LDWT}. The latter is 
related to the proof of \eqref{w-variance} in Corollary \ref{CLTW}.
We emphasize that Theorem \ref{SA} does not
help because almost sure convergence does not say anything for $L^p$ convergence
unless we have bounded random variables, which is definitively not the case here.

It is easy to get that $\mathcal{W}_n'(0)=\frac{1}{n}\int \log\w_n\ d\mu\!\times\!\mu$.
Since $q\mapsto \mathcal{W}_n(q)$ is a convex and, at least, continuously differentiable
function on $(-1,+\infty)$, $\mathcal{W}'(0)=\limn\mathcal{W}_n'(0)$.
Let us show that $\limn\mathcal{W}_n'(0)= h_\mu(T)$.

First observe that convergence in mean of $(\log\w_n)_n$ is equivalent to
showing
$$
\limn\frac1{n}\int \vert\log(\w_n(x,y) \mu([x_1^n])\vert\ d\mu\!\times\!\mu(x,y)=0\,.
$$
Indeed, 
$$
\frac1{n}\int \vert\log\w_n(x,y)-n h_\mu(T)\vert \ d\mu\!\times\!\mu(x,y)\leq
$$
$$
\frac1{n}\int \vert\log(\w_n(x,y)\mu([x_1^n])\vert\ d\mu\!\times\!\mu(x,y)
+
\frac1{n}\int \vert \log\mu([x_1^n]) + n h_\mu(T)\vert\ d\mu(x)\;.
$$ 
The second term goes to zero by using \eqref{BGI}, $P(\phi)=0$ and \eqref{normalization}
and noting that $\int S_n\phi\ d\mu=n\int \phi\ d\mu$ by
$T$-invariance of $\mu$.
We have
$$
\int \vert\log(\w_n(x,y) \mu([x_1^n])\vert\ d\mu\!\times\!\mu(x,y)=
\sum_{a_1^n} \mu([a_1^n]) \int_0^\infty \mu\{\log(\t_{[a_1^n]}\mu([a_1^n]))>t'\}\ dt'
$$
The change of variable $\log t=t'$ leads to
$$
\sum_{a_1^n} \mu([a_1^n]) \int_{1}^\infty \mu\{\t_{[a_1^n]}>t/\mu([a_1^n])\}\ \frac{dt}{t}\,.
$$
By Key-Lemma \ref{MKL}, this integral is finite and bounded between, say, $C_1$ and $C_2$ (independent
of $n$).
Therefore we get the desired result. 

\bigskip

Now turn to the convergence in quadratic mean of $(\log\w_n /n)_n$. Observe that the following identity
holds:
$$
\frac1{n}\int\left( \log\w_n - n h_\mu(T)\right)^2\ d(\mu\!\times\!\mu)=
$$
$$
\frac1{n}\int\left( \log(\w_n \mu([x_1^n]))\right)^2 d(\mu\!\times\!\mu)\ +
\frac1{n}\int\left( \log\mu([x_1^n]+n h_\mu(T)\right)^2 d\mu\ +
$$
$$
-\frac2{n}\int \log(\w_n \mu([x_1^n]))(\log\mu([x_1^n]) +n h_\mu(T))\ d(\mu\!\times\!\mu)
$$
The second term goes to $\sigma^2$ (see formula \eqref{defvariance}
and use \eqref{BGI}).
Hence the proof is done if we show that the two other terms go to $0$
as $n\to\infty$.
Proceeding as before we get that the integral in the first term equals
$$
2\ \sum_{a_1^n} \mu([a_1^n]) \int_{1}^\infty \mu\{\t_{[a_1^n]}>t/\mu([a_1^n])\}\ \frac{\log t}{t}\ dt\,.
$$
Using again Key-lemma \ref{MKL}, we bound the integral from above and
below uniformly in $n$. 
Now consider the integral in the third term which is equal to
$$
\sum_{a_1^n} \mu([a_1^n])\ (\log\mu([a_1^n])+ nh_\mu(T))\ \int 
\log\left(\t_{[a_1^n]}(y)\mu([a_1^n])\right)\ d\mu(y)\,.
$$
We recognize the same integral as above which we know bounded from
above and below uniformly in $n$. The factor in front of the integral
was also treated above in this section.
The proof is finished.

\bigskip

We could prove the same results for $\r_n$ by using Key-lemma \ref{KL2}
and for $\nor_n$ by using the corresponding result found in \cite{maurer}. 

\section*{Acknowledgments}

We would like to thank M. Abadi for useful conversations as well as a
careful reading of an earlier version which contained a mistake.


\end{document}